\documentclass[preprint,12pt]{elsarticle}

\usepackage{graphicx} 
\usepackage{color}

\usepackage{amsmath}
\usepackage{makecell}
\usepackage{soul}

\title{A novel Newton-Raphson style root finding algorithm}

\author[UPEC]{Komi AGBALENYO*}

\affiliation[UPEC]{organization={Faculté des sciences et technologie},
            addressline={61 Av. du Général de Gaulle}, 
            city={Creteil},
            postcode={94000}, 
            country={France}}
\author[SIDAM]{Vincent CAILLIEZ*}
\affiliation[SIDAM]{organization={Poject AP3C},
            addressline={9 Allée Pierre de Fermat}, 
            city={AUBIERE},
            postcode={63170}, ,
            country={France}}
\author[UPEC]{Jonathan CAILLIEZ}

\ead{jonathan.caillez@u-pec.fr}

\begin{document}

\maketitle
\section*{Abstract}
Many problems in applied mathematics require root finding algorithms. Unfortunately, root finding methods have limitations. 
Firstly, regarding the convergence, there is a trade-off between the size of it's domain and it's rate.

Secondly the numerous evaluations of the function and its derivatives penalize the efficiency of high order methods.

In this article, we present a family of high order methods, that require few functional evaluations ( One for each step plus one for each considered derivative at the start of the method), thus increasing the efficiency of the methods.


\section{Introduction}

Finding the root of a function, that is solving the equation~:
\begin{equation}
g(x)=0
\end{equation}
is a widely studied topic in numerical analysis. The methods that are developed to solve this kind of equations are employed widely in the different fields of science. They are used in decision science \cite{decision_applications_2019}, computer aided design \cite{computer_aided_3_2005}, simulation of physical phenomena \cite{zhang_physics_simulation_1_2016} and more domains.

One of the first and most well known methods that has been developed to solve this class of problem is the Newton-Raphson method \cite{ypma1995historical}. This method has been then refined further in order to increase it's order of convergence, it's domain of convergence, or to increase it's efficiency index \cite{said_solaiman_optimal_2019_king,ARGYROS2015336,BEHL201589,CORDERO20102969,SOLEYMANI2012847}.
The efficiency index of a method has been defined in \cite{kung_optimal_1974} as $n^{\frac{1}{q}}$ where q is the number of functional evaluations at each step and n is the order of convergence of the method. 
Unfortunately, high order of convergence methods require more evaluation of the function and its derivatives  so their efficiency is limited. For example the methods presented in these articles \cite{said_solaiman_optimal_2019_king,argyros_convergence_2015_convergence_fourth_order,behl_new_2016_optimal_eighth_order,chun_comparison_2016_comparison_of_eighth_order_methods,amat_third-order_2007_kantorovich_conditions} have high order of convergence but the numerous evaluations of $f$ and its derivatives reduce their efficiency.
\cite{budzko_modifications_2014} shows that there is a trade-off between the size of the domain of convergence and the method rate of convergence. It shows that for some functions as $atan$, the traditional methods don't converge for an initial guess greater than $1.39$. The tricks used to extend the domain of convergence above reduces the rate of convergence of these methods. 
Most of the methods in the literature consider no prior knowledge about the function g 
at the solution point l where $g(l)=0$.
However, under certain circumstances, it is possible to compute the successive derivatives of the function at l, be it by physical reasoning in the case where the function g is derived from a physical simulation, or for a class of function where the derivative can be computed knowing the value of the function. For instance let $g(x)=\sqrt{x}$ then $g'(x)=\frac{1}{2\sqrt{x}}=\frac{1}{2g(x)}$
The work of this paper is divided as follows : the first part presents the proofs and construction of the methods. In the second part presents numerical tests of the methods and compare them to other methods present in the literature.




\section{Methods and Proof}

\subsection{Proof for the method with the use of the first derivative, (second order of convergence)}

Given an equation \begin{equation}
    g(x) = 0
    \label{eq1_j}
\end{equation} 
We note $l$ one of its solutions, we consider g to have a Taylor expansion at $l$ with a non 0 convergence radius.
let  
\begin{equation}
    f(x) = g(x) + x
\end{equation} 
 then $l$ is a fixed point of $f$.
\begin{equation}
    f(l) = g(l) + l=l
\end{equation}
The proposed method starts with the following sequence, while trying to improve its order of convergence :
\begin{equation}
    U_{n+1} = f(U_n)
\end{equation}
Let $ U_n =  l(1 + \xi)  $, then 
\begin{equation} 
    U_{n+1} = f(U_n) = f(l(1 + \xi))
    \label{eq2_j}
\end{equation}.

Using the Taylor series of  $f(x)$ near $l$,  the equation \ref{eq2_j}  becomes 

   \begin{equation}
       U_{n+1}  = f(l) + l\xi f'(l) + \frac{l^2\xi^2}{2}f''(l) + o(l^2\xi^2)  
   \end{equation}

    \begin{equation}
        U_{n+1} - l =  l\xi f'(l) + \frac{l^2\xi^2}{2}f''(l) + o(l^2\xi^2) 
    \end{equation}

In order to accelerate the convergence, we will reduce to zero the first order error $l\xi f'(l) $. To do so, a corrective term is added, and thus the sequence becomes : 

    \begin{equation}
      U_{n+1} = f(U_n) + \alpha(f(U_n) -U_{n})  
    \end{equation}
Which gives the following result using the Taylor series of f
    \begin{equation}
        U_{n+1} = f(l) + l\xi \left[ f'(l) + \alpha(f'(l) - 1)\right]+ o(l\xi)
        \label{eq3_j1}
    \end{equation}
In order to cancel out the $l\xi$ term (thus making the method of order 2), we get :
    \begin{equation}
        \alpha = \frac{f'(l)}{1 - f'(l)}         
    \end{equation}
This gives the following result :

\begin{equation}
    U_{n+1} = f(U_n) + [f(U_n) -U_n]\frac{f'(l)}{1 - f'(l)}
\end{equation}

We replace $f(x)$ by $g(x) + x$ and $f'(x)$ by $g'(x) + 1$, then the sequence to approximate the roots of g(x) can then be expressed as :


\begin{equation}
    U_{n+1} = U_n - \frac{g(U_n)}{g'(l)}
\end{equation}

This method is in theory of order 2, while only requiring one functional evaluation at each step ($g'(l)$ only has to be computed once at the start), it has an efficiency index of 2.

\subsection{Proof for the method with the use of the first and second derivatives, (third order of convergence)}
Given an equation \begin{equation}
    g(x) = 0
\end{equation} 

We note $l$ one of its solutions, we consider $g$ to have a Taylor expansions at $l$ with a non 0 radius of convergence.
Let  
\begin{equation}
    f(x) = g(x) + x
\end{equation} 
$l$ is a fixed point of $f$.
\begin{equation}
    f(l) = g(l) + l=l
\end{equation}
The proposed third order method starts with the following sequence, while trying to improve its order of convergence :
\begin{equation}
    U_{n+1} = f(U_n)
\end{equation}
Let $ U_n =  l(1 + \xi)  $, then 
\begin{equation} 
    U_{n+1} = f(U_n) = f(l(1 + \xi))
    \label{eq3_j2}
\end{equation}.

Using the Taylor series of  $f(x)$ near $l$,  the equation \ref{eq3_j2}  becomes 

\begin{equation}
  U_{n+1}-U_n = l\xi (f'(l) - 1) + \frac{l^2\xi^2}{2}f''(l) + o(l^2\xi^2)  
\end{equation}

\begin{equation}
(U_{n+1}-U_n)^2  = l^2\xi^2 (f'(l) - 1)^2  + o(l^2\xi^2)    
\end{equation}

In order to accelerate the convergence, we will reduce to zero the first order error $l\xi f'(l) $ and the second order error $l^2\xi^2 f''(l)$. To do so, the method becomes:

\begin{equation}
        U_{n+1} = f(U_n) + \alpha(f(U_n) - U_n) + \beta(f(U_n) - U_n)^2 
\end{equation}

\begin{equation}
    U_{n+1} = f(l) + l\xi[f'(l) + \alpha(f'(l) -1)] + 
   \frac{l^2\xi^2}{2}[f''(l)(1 + \alpha) +
   2\beta(f'(l) - 1)^2 ] + 
   o(l^2\xi^2)  
\end{equation}

One deduces these values for $\alpha$ and $\beta$ in order to reduce to zero the $l\xi$ and the $l^2\xi^2$ terms.

\begin{equation}
\begin{array}{cc}
     \alpha = \frac{f'(l)}{1-f'(l)}&  \beta = \frac{-f''(l)}{2(1- f'(l) )^3}\\
\end{array}
\end{equation}

The final sequence is consequently 

\begin{equation}
        U_{n+1} = f(U_n) +[f(U_n) -U_n]\frac{f'(l)}{1-f'(l)}+
    [f(U_n) -U_n]^2\frac{-f''(l)}{2(1 - f'(l))^3}
\end{equation}

By replacing $f(x)$ by $g(x) + x$, $f'(x)$ by $g'(x) + 1$ and  $f''(x)$ by $g''(x)$ we have this sequence

\begin{equation}
    U_{n+1} =  U_n - \frac{g(U_n)}{g'(l)} + \frac{g(U_n)^2g''(l)}{2g'(l)^3}  
\end{equation}

This method is in theory of order 3, while only requiring one functional evaluation at each step ($g'(l)$ and $g''(l)$ only have to be computed once at the start), it has an efficiency index of 3.

\subsection{Proof for the method with the use of the first and second and third derivatives, (fourth order of convergence)}

Given an equation \begin{equation}
    g(x) = 0
    \label{eq1_j2}
\end{equation} 
We note $l$ one of its solutions, we consider g to have a Taylor expansion at $l$ with a non 0 radius of convergence..
let  
\begin{equation}
    f(x) = g(x) + x
\end{equation} 
$l$ is a fixed point of $f$.
\begin{equation}
    f(l) = g(l) + l=l
\end{equation}.
The proposed method starts with the following sequence, while trying to improve its order of convergence:
\begin{equation}
    U_{n+1} = f(U_n)
\end{equation}
Let $ U_n =  l(1 + \xi)  $, then 
\begin{equation} 
    U_{n+1} = f(U_n) = f(l(1 + \xi))
    \label{eq4_j}
\end{equation}.

Using the Taylor series of  $f(x)$ near $l$,  the equation \ref{eq4_j}  becomes 

\begin{equation}
    U_{n+1} = f(l) + l\xi f'(l) + \frac{l^2\xi^2}{2}f''(l) + \frac{l^3\xi^3}{6}f'''(l) +o(l^3\xi^3)
\end{equation}

\begin{equation}
    U_{n+1}-U_n = l\xi(f'(l)-1)+ \frac{l^2\xi^2}{2}f''(l) + \frac{l^3\xi^3}{6}f'''(l) + o(l^3\xi^3)    
\end{equation}

\begin{equation}
    (U_{n+1}-U_n)^2 = l^2\xi^2(f'(l)-1)^2 + l^3\xi^3 f''(l)(f'(l)-1) + o(l^3\xi^3)
\end{equation}

\begin{equation}
    (U_{n+1}-U_n)^3 = l^3\xi^3(f'(l) - 1)^3 + o(l^3\xi^3)    
\end{equation}

In order to accelerate the convergence, we will reduce to zero the first, the second and the third order errors $l\xi f'(l) $,  $l^2\xi^2 f''(l) $ and  $l^3\xi^3 f'''(l) $. To do so, the method becomes: 

\begin{equation}
        U_{n+1} =  f(U_{n}) + \alpha(f(U_{n})  - U_{n}) + \beta(f(U_{n})  - U_{n})^2 
    + \gamma(f(U_{n})  - U_{n})^3
\end{equation}

\begin{multline}
  U{n+1} = f(l) + l\xi[f'(l) + \alpha(f'(l) - 1)] + \\ 
    \quad \frac{l^2\xi^2}{2}[f''(l)(1 + \alpha) + 2\beta(f'(l) - 1)^2] + \\
    \frac{l^3\xi^3}{6}[f'''(l)(1 + \alpha) +6\beta f''(l)(f'(l) - 1) + 6\gamma(f'(l)-1)^3]  
\end{multline}

In order to cancel out the first, the second and the third order errors (thus making the method of order 4), we get 
\begin{equation}
    \begin{array}{c}
\alpha = \frac{f'(l)}{1-f'(l)} \\
\\
  \beta = \frac{-f''(l)}{2(1-f'(l))^3} \\
  \\
  \gamma = \frac{f'''(l)(1 - f'(l)) + 3f''(l)^2}{6(1 -f'(l))^5}
    \end{array}
\end{equation}

We replace $f(x)$ by $g(x) + x$ and $f'(x)$ by $g'(x) + 1$, we do so with the other derivatives and we get this sequence.

\begin{equation}
        U_{n+1} =  U_n - \frac{g(U_n)}{g'(l)} + \frac{g(U_n)^2g''(l)}{2g'(l)^3} \\ -  \frac{g(U_n)^3[3g''(l)^2 -g'(l)g'''(l)]}{6g'(l)^5}
\end{equation}

This method is in theory of order 4, while only requiring one functional evaluation at each step ($g'(l)$, $g''(l)$, $g'''(l)$ only have to be computed once at the start), it has an efficiency index of 4, this is higher than the optimal of $4^\frac{1}{3}$ presented in \cite{said_solaiman_optimal_2019_king}.

\subsection{ Generalization of the proof }
The same method can be applied up to any order of convergence desired. 

Given an equation \begin{equation}
    g(x) = 0
    \label{eq1_j3}
\end{equation} 
We note $l$ one of its solutions, we consider g to have a Taylor expansion at $l$ with a non 0 radius of convergence.
let  
\begin{equation}
    f(x) = g(x) + x
\end{equation} 
$l$ is a fixed point of $f$.
\begin{equation}
    f(l) = g(l) + l=l
\end{equation}.
We assume that the the function \(f\) has a n-th order taylor series

Let,

\begin{equation}
    U_{n+1}-U_{n} =  \sum_{k=0}^{\ n-1} \frac{(l\xi)^k}{k!}f^{k}(l) + o((l\xi)^k) - l(1 + \xi)
\end{equation}

We define after the sequence 

\begin{equation}
    V_{n+1} = U_{n+1} + \sum_{k=1}^{\ n-1} C_{k}(U_{n+1}-U_{n} )^{k}  
\end{equation}

were $C_1$, $C_2$,..., $C_k$ are named constant we will seek in order to reduce to zero all the errors of order smaller than n. Once the formula of all these constants are found we can replace $f(l)$ by $g(l) +l $ and calculate the derivatives accordingly to get the final sequence. This final sequence will give an approximation for the solution $l$ of the equation we want to solve.  
This final sequence will have an efficiency of $n$ because at each step we need only to calculate $f(U_n)$ while the order of convergence is $n$.

\section{Results}
 We will here test our method against other methods present in the literature for different functions and different starting points.
 The computation have been realised on windows 10 computer with and Intel(R) Core(TM) i5-4300M CPU.

\vspace{5mm}

\subsection{Comparison of mehtods for$f(x) =  arctan(x) $ and X0 =$-0.9$ }

\begin{tabular}{|c|c|c|c|} 
\hline  
Method & Time $\mu s$   & Steps & converge\\ 
\hline 
Second order &  1 $\mu s$ &  5 & yes \\  
\hline
Third order &  1 $\mu s$ &  5 & yes \\  
\hline
Fourth order &  2 $\mu s$&  4 & yes \\  
\hline
Newton-Raphson \cite{ehiwario_comparative_2014_comp_new_bisec_secant} &  8 $\mu s$ &  6 & yes \\ 
\hline
Newton two-step \cite{amat_third-order_2007_kantorovich_conditions} &  $<$ 1 $\mu s$ &  4 & yes \\  
\hline
Halley \cite{amat_third-order_2007_kantorovich_conditions} &  & 2 & No \\
\hline
Chebyshev \cite{amat_third-order_2007_kantorovich_conditions} & 2 $\mu s$ & 6 & Yes \\
\hline
\makecell{Derivative free\\ Four order \cite{said_solaiman_optimal_2019_king} }&  &  & No \\ 
\hline
\makecell{Derivative free\\ Eight order \cite{said_solaiman_optimal_2019_king} } &  &  & No \\

\hline 
\end{tabular}

\subsection{Comparison of methods for$f(x) =  arctan(x) $ and X0 =$-10^6$}

\begin{tabular}{|c|c|c|c|} 
\hline  
Method & Time $\mu s$   & Steps & converge\\ 
\hline 
Second order method &   50694 $\mu s$ &  636630 & yes \\  
\hline
Third order method &  46893 $\mu s$ &  636630 & yes \\  
\hline
Fourth order method &  59850 $\mu s$&  349327 & yes \\  
\hline
Newton-Raphson \cite{ehiwario_comparative_2014_comp_new_bisec_secant} &  &   & No \\ 
\hline
Newton two-step \cite{amat_third-order_2007_kantorovich_conditions} &   &   & No \\  
\hline
Halley \cite{amat_third-order_2007_kantorovich_conditions} &  &  & No \\
\hline
Chebyshev \cite{amat_third-order_2007_kantorovich_conditions} &  &  & No \\
\hline
\makecell{Derivative free\\ Four order \cite{said_solaiman_optimal_2019_king} } &  &  & No \\ 
\hline
\makecell{Derivative free\\ Eigth order \cite{said_solaiman_optimal_2019_king} } & &  & No \\ 

\hline 
\end{tabular}

\subsection{Results of the methods for $f(x) =sqrt(\left|x\right|)-4$ and X0 = $-10^{-6}$}
\begin{tabular}{|c|c|c|c|} 
\hline  
Method & Time used  & Steps & converge\\ 
\hline 
Second order method &   11 $\mu s$ &  258 & yes \\  
\hline
Third order method &  $<1$ $\mu s$ &  2 & yes \\  
\hline
Fourth order method &  1 $\mu s$&  2 & yes \\  
\hline
Newton-Raphson \cite{ehiwario_comparative_2014_comp_new_bisec_secant} &  30 $\mu s$ &  252 & Yes \\ 
\hline
Newton two-step \cite{amat_third-order_2007_kantorovich_conditions} &    &  & No \\  
\hline
Halley \cite{amat_third-order_2007_kantorovich_conditions} &  &  & No \\
\hline
Chebyshev \cite{amat_third-order_2007_kantorovich_conditions} & 1 $\mu s$ & 2 & Yes \\
\hline
\makecell{Derivative free\\ Four order \cite{said_solaiman_optimal_2019_king} }&  &  & No \\ 
\hline
\makecell{Derivative free\\ Eigth order \cite{said_solaiman_optimal_2019_king} }&  &  & No \\ 

\hline 
\end{tabular}

\section{Conclusion}

We have presented in this article a family of methods with high efficiency. Their advantage is that they require at each step only the evaluation of $g$.
The others components of the formulas are computed only one time. This gives them for an order of convergence $n$ an efficiency of $n$. Moreover these methods can exhibit a wide domain of convergence. 
Because these methods only use only one functional evaluation at each step, are particularly well suited for the problem of root finding where the valuation of  $g$ and it's derivatives is costly, for instance for physical models.

\section{Acknowledgments}
Vincent CAILLIEZ and Komi AGBALENYO contributed equally to this work.
\section{Citations}
\bibliographystyle{elsarticle-num} 
\bibliography{biblio}

\end{document}